\documentclass[11pt,leqno]{article}

\def\dist{\mathop{\rm dist}}

\begin{document}

\title{Potpourri, 11}

\author{Stephen William Semmes	\\
	Rice University		\\
	Houston, Texas}

\date{}

\maketitle

	Let $(M, d(x, y))$ be a metric space.  Thus $M$ is a nonempty
set and $d(x, y)$ is a nonnegative real-valued function defined for
$x, y \in M$ such that $d(x, y) = 0$ if and only if $x = y$, $d(y, x)
= d(x, y)$ for all $x, y \in M$, and
\begin{equation}
	d(x, z) \le d(x, y) + d(y, z)
\end{equation}
for all $x, y, z \in M$.  If
\begin{equation}
	d(x, z) \le \max(d(x, y), d(y, z))
\end{equation}
for all $x, y, z \in M$, then we say that $d(x, y)$ defines an
ultrametric on $M$.

	A real-valued function $f(x)$ on $M$ is said to be Lipschitz
of order $\alpha$ for some positive real number $\alpha$ if there is a
nonnegative real number $C$ such that
\begin{equation}
	|f(x) - f(y)| \le C \, d(x, y)
\end{equation}
for all $x, y \in M$, which is equivalent to
\begin{equation}
	f(x) \le f(y) + C \, d(x, y)
\end{equation}
for all $x, y \in M$.  Notice that this holds with $C = 0$ if and only
if $f$ is constant on $M$.

	If $p \in M$ and $f(x) = d(x, p)$, then $f$ is Lipschitz of
order $1$ with $C = 1$ by the triangle inequality.  More generally, if
$A$ is a nonempty subset of $M$, then
\begin{equation}
	\dist(x, A) = \inf \{d(x, y) : y \in A\}
\end{equation}
is Lipschitz of order $1$ with $C = 1$.

	Every Lipschitz function on the real line of order $\alpha >
1$ is constant.  One way to see this is to observe that such a
function is differentiable at every point with derivative equal to
$0$.  More generally, every Lipschitz function on ${\bf R}^n$ of order
strictly larger than $1$ is constant.

	If $a$, $b$ are nonnegative real numbers and $q$ is a positive
real number, then
\begin{equation}
	\max (a, b) \le (a^q + b^q)^{1/q}
\end{equation}
and hence
\begin{equation}
	(a^r + b^r)^{1/r} \le \max (a, b)^{(r-q)/r} (a^q + b^q)^{1/r}
		\le (a^q + b^q)^{1/q}
\end{equation}
when $r > q$.  In particular, if $0 < q \le 1$, then
\begin{equation}
	(a + b)^q \le a^q + b^q.
\end{equation}

	Let $(M, d(x, y))$ be a metric space and let $q$ be a positive
real number with $q \le 1$.  One can check that $d(x, y)^q$ is also a
metric on $M$.

	Now suppose that $d(x, y)$ is an ultrametric on $M$ and that
$q$ is any positive real number.  One can check $d(x, y)^q$ is an
ultrametric on $M$ too.  In both situations the topology on $M$
determined by $d(x, y)^q$ is the same as the one determined by $d(x,
y)$.

	Suppose that $(M, d(x, y))$ is a metric space and that
\begin{equation}
	d(x, y) = \rho(x, y)^a
\end{equation}
for some metric $\rho(x, y)$ on $M$ and some $a > 0$.  A real-valued
function $f(x)$ on $M$ is Lipschitz of order $\alpha > 0$ on $M$ with
constant $C$ with respect to the metric $d(x, y)$ if and only if
$f(x)$ is Lipschitz of order $a \, \alpha$ on $M$ with constant $C$
with respect to $\rho(x, y)$.

	If $a > 1$, then it follows that there are plenty of Lipschitz
functions on $M$ of order $a$ with respect to $d(x, y)$, since there
are plenty of Lipschitz functions of order $1$ on $M$ with respect to
$\rho(x, y)$.  For example, these functions separate points on $M$,
which is to say that for any pair of distinct points in $M$ there is
such a function which takes different values at the two points.

	As a special case, consider the real line ${\bf R}$ with the
metric $|x - y|^a$ for some $a \in (0, 1)$.  A Lipschitz function on
${\bf R}$ of order $1$ with respect to the standard metric $|x - y|$
is Lipschitz of order $1/a > 1$ with respect to $|x - y|^a$.  For each
positive integer $n$ we can consider ${\bf R}^n$ with the metric $|x -
y|^a$, where $|x - y|$ is the standard metric, and every Lipschitz
function on ${\bf R}^n$ of order $1$ with respect to the standard
metric is a Lipschitz function of order $1/a$ with respect to $|x -
y|^a$.

	Now consider ${\bf R}^2$ where the distance between $x = (x_1,
x_2)$ and $y = (y_1, y_2)$ is defined to be
\begin{equation}
	|x_1 - y_1| + |x_2 - y_2|^{1/2}.
\end{equation}
One can check that this defines a metric on ${\bf R}^2$ which is
compatible with the standard topology.

	This kind of space was considered by my colleague Frank Jones
in connection with the heat equation and related singular integral
operators.  It is also an instance of a ``Rickman's rug'' as discussed
in \cite{Juha3}, with very interesting properties concerning
quasiconformal mappings and geometry.

	As on any metric space, there are plenty of Lipschitz
functions on this space of order $\alpha \le 1$.  There are
nonconstant Lipschitz functions of order $\alpha \in (1, 2]$, which
are constant in the first variable.  Every Lipschitz function on this
space of order strictly larger than $2$ is constant.

	It is convenient to define $\|x\|$ for $x \in {\bf R}^2$ by
\begin{equation}
	\|x\| = |x_1| + |x_2|^{1/2},
\end{equation}
so that the distance function on ${\bf R}^2$ defined above can be
described as $\|x - y\|$.  Notice that
\begin{equation}
	\|x + y\| \le \|x\| + \|y\|
\end{equation}
for all $x, y \in {\bf R}^2$, and if we define dilations $\delta_r$ on
${\bf R}^2$ for $r > 0$ by
\begin{equation}
	\delta_r(x) = (r \, x_1, r^2 \, x_2),
\end{equation}
then
\begin{equation}
	\|\delta_r(x)\| = r \, \|x\|
\end{equation}
for all $x \in {\bf R}^2$ and $r > 0$.

	One can check that the measure of a ball of radius $t$ in
${\bf R}^2$ with respect to this metric and using standard Lebesgue
measure is equal to a positive constant times $t^3$, where the
constant is simply the measure of the ball centered at the origin with
radius $1$.  Also, every open ball in ${\bf R}^2$ with respect to this
metric is homeomorphic to the standard unit disk in ${\bf R}^2$.

	Let $(M, d(x, y))$ be a metric space, and let $\epsilon$,
$\lambda$ be positive real numbers.  A finite sequence of points $z_0,
\ldots, z_l$ in $M$ will be called an $(\epsilon, \lambda)$-chain if
\begin{equation}
	d(z_i, z_{i+1}) < \epsilon \, \lambda
\end{equation}
for all $i$, $0 \le i < l$, and
\begin{equation}
	\sum_{i=0}^{l-1} d(z_i, z_{i+1}) \le \lambda.
\end{equation}
In other words, the size of each step in the chain is less than
$\epsilon$, and the total length of the chain is less than or equal to
$\lambda$.

	Suppose that $f(x)$ is a real-valued function on $M$ which is
Lipschitz of order $\alpha > 1$ with constant $C$.  If $z_0, \ldots,
z_l$ is an $(\epsilon, \lambda)$-chain in $M$, then we have that
\begin{equation}
	|f(z_0) - f(z_l)| \le \sum_{i=0}^{l-1} |f(z_i) - f(z_{i+1})|
  		\le C \, d(z_i, z_{i+1})^\alpha 
			\le C \, \epsilon^{\alpha - 1} \, \lambda^\alpha.
\end{equation}

	On ${\bf R}^n$ with the standard metric, each pair of points
can be connected by an $(\epsilon, \lambda)$-chain where $\epsilon >
0$ is arbitrary and $\lambda$ is equal to the distance between the two
points, and this gives another way to see that Lipschitz functions of
order strictly larger than $1$ must be constant.  In the real line
there are Cantor sets on which there are plenty of locally constant
functions, and hence plenty of functions which are Lipschitz of any
positive order, and where the gaps in the Cantor set are small enough
so that there are a lot of interesting $(\epsilon, \lambda)$-chains
with $\epsilon$ small, which lead to significant restrictions on the
behavior of Lipschitz functions of order strictly larger than $1$.

	Let $(M_1, d_1(x, y))$ and $(M_2, d_2(u, v))$ be metric
spaces.  A mapping $\phi : M_1 \to M_2$ is said to be bilipschitz with
constant $C \ge 1$ if
\begin{equation}
	C^{-1} \, d_1(x, y) \le d_2(\phi(x), \phi(y)) \le C \, d_1(x, y)
\end{equation}
for all $x, y \in M_1$.  If there is a bilipschitz mapping of $M_1$
onto $M_2$, then we say that $M_1$ and $M_2$ are bilipschitz
equivalent.

	Notice that $C = 1$ in this condition if and only if $\phi :
M_1 \to M_2$ is an isometry, which is to say that
\begin{equation}
	d_2(\phi(x), \phi(y)) = d_1(x, y)
\end{equation}
for all $x, y \in M_1$.  If there is an isometry of $M_1$ onto $M_2$,
then we say that $M_1$, $M_2$ are isometrically equivalent.

	Let $(M, d(x, y))$ be a metric space, and let $E$ be a subset
of $M$.  We say that $E$ is porous in $M$ if there is a $C > 0$ so
that for each $x \in M$ there is a $y \in M$ with $d(x, y) \le r$ and
$d(y, z) \ge C^{-1} \, r$ for all $z \in E$.

	Clearly a porous set has empty interior.  Observe that the
closure of a porous set is again porous, and with the same constant
$C$.  Thus a porous set is nowhere dense.

	Let us take $M$ to be ${\bf R}^n$ with the standard metric.
It is easy to see that a porous set $E$ has measure $0$, and indeed
even if one takes the closure of $E$ there are no points of density.

	We can strengthen this assertion as follows.  If $Q$ is a cube
in ${\bf R}^n$ and $L \ge 2$ is an integer, then we can express $Q$ as
the union of $L^n$ cubes with sidelength equal to $1/L$ times the
sidelength in the usual way, by subdividing the sides of $Q$ into $L$
segments of equal length.  A subset $E$ of ${\bf R}^n$ is porous if
and only if there is an $L \ge 2$ such that for each cube $Q$ in ${\bf
R}^n$ there is a subcube $Q_1$ of $Q$ which is disjoint from $E$ and
which is one of the $L^n$ subcubes with sidelength equal to $1/L$
times the sidelength of $Q$ just described.

	Using this one can check that for each cube $Q$ in ${\bf R}^n$
and each positive integer $k$ one can cover $E \cap Q$ by $\le (L^n -
1)^k$ subcubes of $Q$ with sidelength equal to $L^{-k}$ times the
sidelength of $Q$.  It follows that the Hausdorff dimension of $E$ is
less than or equal to $(\log (L^n - 1)) / (\log L) < n$.


\begin{thebibliography}{53}

\bibitem {A-V1} P.~Alestalo and J.~V\"ais\"al\"a, {\it Quasisymmetric
embeddings of products of cells into the Euclidean Space}, Annales
Academi\ae Scientiarum Fennic\ae Series A I Mathematica {\bf 19}
(1994), 375--392.

\bibitem {A-T} L.~Ambrosio and P.~Tilli, {\it Topics on Analysis in
Metric Spaces}, Oxford University Press, 2004.

\bibitem {A1} P.~Assouad, {\it Espaces M\'etriques, Plongements,
Facteurs}, Th\`ese de Doctorat (January, 1977), Universit\'e de Paris
XI, 91405 Orsay, France.

\bibitem {A2} P.~Assouad, {\it \'Etude d'une dimension m\'etrique
li\'ee \`a la possibilit\'e de plongement dans ${\bf R}^n$}, Comptes
Rendus de l'Acad\'emie des Sciences Paris, S\'er.~A {\bf 288} (1979),
731--734.

\bibitem {A3} P.~Assouad, {\it Plongements Lipschitziens dans ${\bf
R}^n$}, Bulletin de la Soci\'et\'e Math\'ematique de France {\bf 111}
(1983), 429--448.

\bibitem {C-W-1} R.~Coifman and G.~Weiss, {\it Analyse Harmonique
Non-Commutative sur Certains Espaces Homog\`enes}, Lecture Notes in
Mathematics {\bf 242}, Springer-Verlag, 1971.

\bibitem {C-W-2} R.~Coifman and G.~Weiss, {\it Extensions of Hardy
Spaces and their Use in Analysis}, Bulletin of the American
Mathematical Society {\bf 83} (1977), 569--645.

\bibitem {F-S} G.~Folland and E.~Stein, {\it Hardy Spaces on
Homogeneous Groups}, Princeton University Press, 1982.

\bibitem {Gol} R.~Goldberg, {\it Methods of Real Analysis}, 2nd
edition, Macmillan, 1976.

\bibitem {Juha1} J.~Heinonen, {\it Calculus on Carnot Groups}, in {\it
Fall School in Analysis}, Jyv\"askyl\"a, 1994, Reports of the
Department of Mathematics and Statistics, University of Jyv\"akyl\"a
{\bf 68}, 1--31, 1995.

\bibitem {Juha2} J.~Heinonen, {\it Lectures on Analysis on Metric
Spaces}, Springer-Verlag, 2001.

\bibitem {Juha3} J.~Heinonen, {\it Geometric Embeddings of Metric
Spaces}, Reports of the Department of Mathematics and Statistics,
University of Jyv\"askyl\"a {\bf 90}, 2003.

\bibitem {H-S} E.~Hewitt and K.~Stromberg, {\it Real and Abstract
Analysis}, Springer-Verlag, 1975.

\bibitem {H-W} W.~Hurewicz and H.~Wallman, {\it Dimension Theory},
Princeton University Press, 1941.

\bibitem {J1} F.~Jones, {\it A class of singular integrals}, American
Journal of Mathematics {\bf 86} (1964), 441--462.

\bibitem {J2} F.~Jones, {\it Lebesgue Integration on Euclidean Space},
Jones and Bartlett Publishers, 1993.

\bibitem {Jean-Lin} J.-L.~Journ\'e, {\it Calder\'on--Zygmund
Operators, Pseudodifferential Operators and the Cauchy Integral of
Calder\'on}, Lecture Notes in Mathematics {\bf 994}, Springer-Verlag,
1983.

\bibitem {K} J.~Kigami, {\it Analysis on Fractals}, Cambridge
University Press, 2001.

\bibitem {K1} S.~Krantz, {\it Real Analysis and Foundations},
CRC Press, 1991.

\bibitem {K2} S.~Krantz, {\it A Panorama of Harmonic Analysis},
Mathematical Association of America, 1999.

\bibitem {K3} S.~Krantz, {\it Function Theory of Several Complex
Variables}, AMS Chelsea Publishing, 2001.

\bibitem {K4} S.~Krantz, {\it Complex Analysis: The Geometric
Viewpoint}, 2nd edition, Mathematical Association of America, 2004.

\bibitem {Lu} J.~Luukkainen, {\it Assouad dimension: Antifractal
metrization, porous sets, and homogeneous measures}, Journal of the
Korean Mathematical Society {\bf 35} (1998), 23--76.

\bibitem {Lu-ML} J.~Luukkainen and H.~Movahedi-Lankarani, {\it Minimal
bi-Lipschitz embedding dimension of ultrametric spaces}, Fundamenta
Mathematicae {\bf 144}, 181--193.

\bibitem {M-S} R.~Mac\'{\i}as and C.~Segovia, {\it Lipschitz functions
on spaces of homogeneous type}, Advances in Mathematics {\bf 33}
(1979), 257--270.

\bibitem {Mat} P.~Mattila, {\it Geometry of Sets and Measures in
Euclidean Spaces: Fractals and Rectifiability}, Cambridge University
Press, 1995.

\bibitem {Res} Y.~Reshetnyak, {\it Space Mappings with Bounded
Distortion}, American Mathematical Society, 1989.

\bibitem {Ric} S.~Rickman, {\it Quasiregular Mappings},
Springer-Verlag, 1993.

\bibitem {Roy} H.~Royden, {\it Real Analysis}, 3rd edition, Macmillan,
1988.

\bibitem {Ru1} W.~Rudin, {\it Principles of Mathematical Analysis},
3rd edition, McGraw-Hill, 1976.

\bibitem {Ru2} W.~Rudin, {\it Function Theory in the Unit Ball of
${\bf C}^n$}, Springer-Verlag, 1980.

\bibitem {St1} E.~Stein, {\it Singular Integrals and Differentiability
Properties of Functions}, Princeton University Press, 1970.

\bibitem {St2} E.~Stein, {\it Harmonic Analysis: Real-Variable
Methods, Orthogonality, and Oscillatory Integrals}, with the
assistance of T.~Murphy, Princeton University Press, 1993.

\bibitem {S-W} S.~Staples and L.~Ward, {\it Quasisymmetrically thick
sets}, Annales Academi\ae Scientiarum Fennic\ae Mathematic\ae {\bf 23}
(1998), 151--168.

\bibitem {TrV} D.~Trotsenko and J.~V\"ais\"al\"a, {\it Upper sets and
quasisymmetric maps}, Annales Academi\ae Scientiarum Fennic\ae
Mathematic\ae {\bf 24} (1999), 465--488.

\bibitem {Tu} P.~Tukia, {\it A quasiconformal group not isomorphic to
a M\"obius group}, Annales Academi\ae Scientiarum Fennic\ae Series A I
Mathematica {\bf 6} (1981), 149--160.

\bibitem {TuV} P.~Tukia and J.~V\"ais\"al\"a, {\it Quasisymmetric
embeddings of metric spaces}, Annales Academi\ae Scientiarum Fennic\ae
Mathematica Series A I Mathematica {\bf 5} (1980), 97--114.

\bibitem {V1} J.~V\"ais\"al\"a, {\it Lectures on $n$-Dimensional
Quasiconformal Mappings}, Lecture Notes in Mathematics {\bf 229},
Springer-Verlag, 1971.

\bibitem {V2} J.~V\"ais\"al\"a, {\it Quasisymmetric embeddings in
Euclidean spaces}, Transactions of the American Mathematical Society
{\bf 264} (1981), 191--204.

\bibitem {V3} J.~V\"ais\"al\"a, {\it Porous sets and quasisymmetric
maps}, Transactions of the American Mathematical Society {\bf 299}
(1987), 525--533.

\bibitem {V4} J.~V\"ais\"al\"a, {\it Quasisymmetric maps of products
of curves into the plane}, Revue Roumaine de Math\'ematiques Pures et
Appliqu\'ees {\bf 33} (1988), 147--156.

\bibitem {V5} J.~V\"ais\"al\"a, {\it Quasisymmetric maps}, in {\it
Holomorphic Functions and Moduli}, Volume 1, 187--192, Mathematical
Sciences Research Institute Publications {\bf 10}, Springer-Verlag,
1988.

\bibitem {V6} J.~V\"ais\"al\"a, {\it Free quasiconformality in Banach
Spaces I}, Annales Academi\ae Scientiarum Fennic\ae Series A I
Mathematica {\bf 15} (1990), 355--379.

\bibitem {V7} J.~V\"ais\"al\"a, {\it Free quasiconformality in Banach
Spaces II}, Annales Academi\ae Scientiarum Fennic\ae Series A I
Mathematica {\bf 16} (1991), 255--310.

\bibitem {V8} J.~V\"ais\"al\"a, {\it Free Quasiconformality in Banach
Spaces III}, Annales Academi\ae Scientiarum Fennicae Series A I
Mathematica {\bf 17} (1992), 393--408.

\bibitem {V9} J.~V\"ais\"al\"a, {\it Free quasiconformality in Banach
Spaces IV}, in {\it Analysis and Topology}, 697--717, World Scientific
Publishing, 1998.

\bibitem {V10} J.~V\"ais\"al\"a, {\it The free quasiworld: Freely
quasiconformal and related maps in Banach Spaces}, in {\it
Quasiconformal Geometry and Dynamics}, 55--118, Banach Center
Publications {\bf 48}, 1999.

\bibitem {VVW} J.~V\"ais\"al\"a, M.~Vourinen, and H.~Wallin, {\it
Thick sets and quasisymmetric maps}, Nagoya Mathematical Journal {\bf
135} (1994), 121--148.

\bibitem {W} P.~Wojtaszczyk, {\it Banach Spaces for Analysts},
Cambridge University Press, 1991.

\bibitem {Z} A.~Zygmund, {\it Trigonometric Series}, 3rd edition,
Cambridge University Press, 2002.


\end{thebibliography}
\end{document}